\documentclass{amsart}

\newtheorem{theorem}{Theorem}[section]
\newtheorem{propo}[theorem]{Proposition}
\newtheorem{lemma}[theorem]{Lemma}
\newtheorem{mainthm}{Main Theorem}

\theoremstyle{definition}

\theoremstyle{remark}

\numberwithin{equation}{section}

\begin{document}

\title[Arithmeticity of Lie foliations]{Arithmeticity of totally geodesic Lie foliations with locally symmetric leaves}

\author{R. Quiroga-Barranco}

\address{Centro de Investigaci\'on en Matem\'aticas,
A.~P.~402, Guanajuato, Gto., C.P.~36000, M\'exico}
\email{quiroga@cimat.mx}

\thanks{Supported by SNI-M\'exico and Conacyt Grant 44620}



\begin{abstract}
Zimmer \cite{Zimmer-Lie} proved that, on a compact manifold, a
foliation with a dense leaf, a suitable leafwise Riemannian
symmetric metric and a transverse Lie structure has arithmetic
holonomy group. In this work we improve such result for totally
geodesic foliations by showing that the manifold itself is
arithmetic. This also gives a positive answer, for some special
cases, to a conjecture of E. Ghys \cite{Molino}.
\end{abstract}

\maketitle

\section{Introduction}
In this work we want to consider a compact manifold $M$ with a
smooth foliation $\mathcal{F}$. It is well known that for such
manifolds, the introduction of transverse structures has played a
fundamental role in the theory of foliations. In particular, Lie
foliations, i.e. transversely modeled on a Lie group (see Section
\ref{sec-prel} for more details), are very useful when studying all
sorts of transverse geometries. A remarkable case is given by the
work of Molino on Riemannian foliations \cite{Molino}. A natural
problem to consider is the classification of Lie foliations on
compact manifolds. Nevertheless, this seems to be a problem too
general to solve without further considerations.

The following construction provides a fundamental example of a Lie
foliation. Let $\pi : L \rightarrow H$ be a surjective homomorphism
between connected Lie groups. Then the fibers of $\pi$ yield a
foliation on $L$ which is a Lie foliation modeled on $H$. If
$\Gamma$ is a cocompact discrete subgroup of $L$, then its left
action preserves the Lie foliation on $L$ and so there is an induced
Lie foliation on the compact manifold $\Gamma\backslash L$.
Moreover, if $K$ is a compact subgroup of the kernel of $\pi$, then
the double coset $\Gamma\backslash L/K$ is a compact manifold
carrying a Lie foliation modeled in $H$ as well. We will refer to
these double coset examples as homogeneous Lie foliations. E. Ghys
has conjectured (see Appendix~E to \cite{Molino}) that every Lie
foliation can be obtained from a fibration over a homogeneous Lie
foliation, as long as the model group $H$ has no compact factors; we
refer to page 301 of \cite{Molino} for the precise statement.

On the other hand, one can study Lie foliations that also carry some
sort of leafwise structure. Zimmer \cite{Zimmer-Lie} did so by
considering a compact foliated manifold carrying both a transverse
Lie structure and a leafwise Riemannian metric so that the leaves
are covered by a symmetric space of noncompact type and higher rank.
It was established in \cite{Zimmer-Lie} the following arithmeticity
for the holonomy group of such foliations. We refer to
\cite{Zimmer-Lie} for the complete statement. We observe that the
result in \cite{Zimmer-Lie} is stated for foliations with a simply
connected dense leaf, but the existence of one such leaf is
guaranteed just from the presence of a dense one by the main result
in \cite{Stuck-Zimmer}.

\begin{theorem}[Zimmer \cite{Zimmer-Lie}]
\label{Zimmer-hol-thm} Let $H$ be a connected Lie group and
$\mathcal{F}$ an $H$-foliation (i.e. with a transverse Lie structure
modeled on $H$) of a compact manifold $M$. Suppose that
$\mathcal{F}$ has a dense leaf. Assume there is a leafwise
Riemannian metric on $M$ such that each leaf is isometrically
covered by an irreducible symmetric space of noncompact type and
rank at least $2$. Then $H$ is semisimple and the holonomy group
$\Lambda\subset H$ is a dense arithmetic subgroup, in other words,
$\operatorname{Ad}_H(\Lambda)$ is commensurable with the image in
$\operatorname{Ad}_H(H)$ of an arithmetic lattice $\Gamma$ of a
semisimple Lie group $L$ with respect to a smooth surjection $L
\rightarrow \operatorname{Ad}_H(H)$.
\end{theorem}

Once such notions of arithmeticity had been considered, we can
specialize the above construction of a homogeneous Lie foliation as
follows. Suppose both $L$ and $H$ are semisimple Lie groups with
finite center, and that $\Gamma$ is an arithmetic irreducible
lattice. Then the double coset $\Gamma\backslash L/K$ is called an
arithmetic manifold, and so we will say that the foliation
constructed above from the epimorphism $\pi$ is an arithmetic
homogeneous Lie foliation. Note that if $K$ is a maximal compact
subgroup of the kernel of $\pi$, say $G$, then $\Gamma\backslash
L/K$ carries a leafwise Riemannian metric so that each leaf is
isometrically covered by the symmetric space $G/K$.

Given Ghys' conjecture and Zimmer's arithmeticity result, one can
consider the problem of determining the relation between Lie
foliations with locally symmetric leaves as in Theorem
\ref{Zimmer-hol-thm} and arithmetic homogeneous Lie foliations. The
main goal of this work is to prove that any foliation as in Theorem
\ref{Zimmer-hol-thm} which is also totally geodesic is, up to a
finite covering, an arithmetic homogeneous Lie foliation. From now
on, $G$ denotes a connected noncompact simple Lie group with trivial
center and $X_G = G/K$ the symmetric space associated to $G$, where
$K$ is a maximal compact subgroup of $G$.

\begin{mainthm}
Let $M$ be a compact manifold carrying a smooth foliation
$\mathcal{F}$ that admits a transverse Lie structure and a
Riemannian metric for which $\mathcal{F}$ is totally geodesic and
such that every leaf is isometrically covered by the symmetric space
$X_G$. Suppose that $\mathcal{F}$ has a dense leaf and $X_G$ is
irreducible with rank at least $2$. Then, up to a finite covering,
$M$ is an arithmetic manifold and $\mathcal{F}$ is an arithmetic
homogeneous Lie foliation. More precisely, there exist:
\begin{enumerate}
\item a finite covering map $\pi : \widehat{M} \rightarrow M$,
\item a connected semisimple Lie group $H$ with finite center,
\item an arithmetic irreducible lattice $\Gamma$ of $H \times G$,
and
\item a diffeomorphism $\varphi : \Gamma\backslash
(H\times G)/K = \Gamma\backslash (H\times G/K) \rightarrow
\widehat{M}$, where $K$ is a maximal compact subgroup of $G$,
\end{enumerate}
such that $\pi\circ\varphi$ maps the foliation of $\Gamma\backslash
(H\times G/K)$ induced by the factor $G/K$ of $H\times G/K$
diffeomorphically onto the foliation of $M$, and so that it
preserves both the leafwise Riemannian symmetric metric and the
transverse Lie structure.
\end{mainthm}

This result provides a positive answer to Ghys' conjecture in a
special case. It also relates to Zimmer's arithmeticity theorem for
the holonomy of a Lie foliation since it shows that not just the
holonomy but the foliation and the manifold are arithmetic if we
consider the case of totally geodesic foliations.

The proof of this result use our previous work on actions of
pseudoRiemannian manifolds (see \cite{Quiroga-Annals}) and build on
Theorem \ref{Zimmer-hol-thm}. It is already well known that for a
manifold $M$ as in the Main Theorem, there is a manifold $M^*$ acted
upon by $G$ and that fibers, as a principal bundle with compact
group, over $M$; this is in fact one of the main objects considered
in \cite{Zimmer-Lie}. Here we show that on such a manifold $M^*$ we
can build a $G$-invariant pseudoRiemannian metric thus allowing us
to apply the results from \cite{Quiroga-Annals}. As it will be clear
from the proof, the compactness in the Main Theorem is only used to
ensure completeness of the transverse structures considered. Hence,
our Main Theorem holds for $M$ noncompact carrying a complete
transverse Lie structure and a holonomy invariant transverse smooth
measure.

\section{Preliminaries on transverse and tangential structures}
\label{sec-prel} In this section we assume that $M$ is a manifold
equipped with a smooth foliation $\mathcal{F}$. Also, we denote with
$T\mathcal{F}$ the tangent bundle to the leaves.

We start by describing some basic facts about transverse structures
on $\mathcal{F}$. For these, the fundamental notion is that of the
transverse frame bundle $L_T(\mathcal{F})$ of the foliation which is
defined as the linear frame bundle of $TM / T\mathcal{F}$.

A foliate vector field on $M$ is a smooth vector field $X$ on $M$
such that, for every smooth vector field $Y$ over $M$ tangent to the
leaves, the Lie bracket $[X,Y]$ is tangent to the leaves as well. A
transverse field for the foliated manifold $M$ is a section
$\bar{X}$ of $TM / T\mathcal{F}$ that can be seen as the image under
the natural projection $TM \rightarrow TM / T\mathcal{F}$ of a
foliate vector field. The vector space of transverse vector fields,
denoted with $\mathfrak{l}(M,\mathcal{F})$, is easily seen to
inherit a Lie algebra structure from the Lie brackets of foliate
vector fields (see \cite{Molino} for details). We define a
transverse parallelism of the foliation $\mathcal{F}$ as a family of
pointwise linearly independent transverse fields $\bar{X}_1,\dots,
\bar{X}_q$, where $q$ is the codimension of $\mathcal{F}$ in $M$.
Clearly, a transverse parallelism defines a trivialization of $TM /
T\mathcal{F}$ and thus a reduction of $L_T(\mathcal{F})$ to the
identity group. Conversely, a reduction of $L_T(\mathcal{F})$ to the
identity group defines a transverse parallelism if the total space
of such reduction is saturated for the lifted foliation on
$L_T(\mathcal{F})$. For the definition of the lifted foliation, the
last claim and basic associated notions we refer to \cite{Molino}. A
transverse parallelism is called complete if the transverse vector
fields that define it are the image of complete foliate vector
fields.

From these remarks, a transverse parallelism defines reductions of
$L_T(\mathcal{F})$ to both $O(q)$ and
$\operatorname{SL}^\pm(q;\mathbb{R})$ which are saturated with
respect to the lifted foliation. Hence, from a transverse
parallelism, we obtain a transverse Riemannian structure and a
holonomy invariant transverse smooth measure. Note that the former
is given by declaring the transverse fields that provide the
parallelism as pointwise orthonormal. Furthermore, the transverse
measure defined by the transverse Riemannian metric and the
transverse measure given by the
$\operatorname{SL}^\pm(q;\mathbb{R})$-reduction of
$L_T(\mathcal{F})$ are easily seen to be the same.

From now on, for a foliation $\mathcal{F}$, a transverse Lie
structure modeled on a connected Lie group $H$ (or an $H$-foliation)
will be given by an open cover $\{ U_i \}_i$ of $M$ and smooth
submersions $f_i : U_i \rightarrow H$, compatible with the
foliation, onto open subsets of $H$ so that for every $i,j$ there is
some $x_{ij} \in H$ for which we have:
$$
f_i = L_{x_{ij}}\circ f_j
$$
on $U_i\cap U_j$, where $L_x$ is the left translation in $H$ by a
given $x\in H$.

Alternatively, an $H$-foliation on $\mathcal{F}$ can be defined as a
transverse parallelism $\bar{X}_1,\dots, \bar{X}_q$ whose linear
span over $\mathbb{R}$ is a Lie algebra isomorphic to the Lie
algebra $\mathfrak{h}$ of $H$ (see \cite{Molino}). Given this
remark, we say that a transverse Lie structure is complete if the
completeness condition holds for the corresponding transverse
parallelism. Also, observe that the transverse Riemannian structure
given by the parallelism of a transverse Lie structure modeled on
the group $H$ corresponds to a choice of a left invariant metric on
$H$.

Finally, we refer to \cite{Molino, Zimmer-Lie} and their references
for a more precise account of the basic notions associated to
transverse Lie structures, including those of a development,
holonomy representation and holonomy group of a Lie foliation. We
will use freely the basic facts of the theory of transverse Lie
structures, providing references when necessary.

For the notion of tangential structure for totally geodesic
foliations, dual to the notion of transverse structure, the
fundamental object is $L_{tg}(\mathcal{F})$, the linear frame bundle
of the vector bundle $T\mathcal{F}$. For the rest of this section we
assume that $M$ carries a Riemannian metric for which $\mathcal{F}$
is totally geodesic but with otherwise arbitrary geometry along the
leaves.

One of the most useful properties of totally geodesic foliations is
the fact that the leaves are isometrically unchanged as we move
perpendicularly to them. More precisely, if $\gamma : [0,1]
\rightarrow M$ is a curve starting at $x$ and everywhere
perpendicular to $\mathcal{F}$, then there exists a unique family of
isometries $(\psi_t)_{t\in [0,1]}$ such that $\psi_t : V_x
\rightarrow V_{\gamma(t)}$, where $V_{\gamma(t)}$ is a neighborhood
of $\gamma(t)$ in the leaf that contains $\gamma(t)$, $\psi_0$ is
the identity map and the curves $y \mapsto \psi_t(y)$ are everywhere
perpendicular to $\mathcal{F}$ for every $y\in V_x$ (see
\cite{Johnson-Whitt} for details). The isometries $(\psi_t)_t$ are
called the elements of horizontal holonomy associated to $\gamma$.

Denote with $T\mathcal{F}^\perp$ the orthogonal complement of
$T\mathcal{F}$. Then, the vector bundle $T\mathcal{F}^\perp$ is
lifted to a vector subbundle $\mathcal{H}$ of
$TL_{tg}(\mathcal{F})$, called the lifted horizontal bundle, as
follows (see \cite{Blumenthal-Hebda2,Cairns1} for further details).
For a given $v \in T_x\mathcal{F}^\perp$ choose a curve $\gamma :
[0,1] \rightarrow M$ starting at $x$ with velocity $v$ and
everywhere perpendicular to $\mathcal{F}$. Let $(\psi_t)_t$ be the
elements of horizontal holonomy associated to $\gamma$. If $\alpha
\in L_{tg}(\mathcal{F})$ is an element in the fiber over $x$, then
the curve $\overline{\gamma}(t)=d(\psi_t)_x\circ\alpha$ lies in
$L_{tg}(\mathcal{F})$ and its velocity vector $\overline{v} =
\overline{\gamma}'(0)$ is called the horizontal lift of $v$ at
$\alpha$. The set of all such horizontal lifts of elements of
$T_x\mathcal{F}^\perp$ to the point $\alpha$ is by definition the
fiber $\mathcal{H}_\alpha$ of the bundle $\mathcal{H}$. The lifted
horizontal bundle is used to define tangential structures. More
precisely, for a subgroup $F$ of $\operatorname{GL}(p)$ (where $p$
is the dimension of the leaves of $\mathcal{F}$), a tangential
$F$-structure for $(M,\mathcal{F})$ is an $F$-reduction $Q$ of
$L_{tg}(\mathcal{F})$ such that $\mathcal{H}_\alpha \subset T_\alpha
Q$ for every $\alpha \in Q$. For $F=\{e\}$, the tangential structure
is called a tangential parallelism. Another example is given by
$O_{tg}(\mathcal{F})$, the principal fiber bundle of orthonormal
frames for $T\mathcal{F}$, which defines a tangential
$O(p)$-structure. On the other hand, the pull-back in
$O_{tg}(\mathcal{F})$ of the leaves of $\mathcal{F}$ defines a
foliation $\mathcal{F}_{tg}$. The bundle $T\mathcal{F}_{tg}$ can be
endowed with a natural parallelism given by the standard horizontal
vector fields with respect to the Levi-Civita connection along the
leaves and the vector fields defined by the $O(p)$-action on
$O_{tg}(\mathcal{F})$.

\section{A $G$-space associated to $M$}
\label{sec-M^*} In this section we recall the construction of a
manifold fibered over $M$ by a compact group and on which there is
an action of the group $G$, for $M$ and $G$ as in the Main Theorem.
Our main contribution here is the construction of a pseudoRiemannian
metric which turns out to be $G$-invariant. The next result is
essentially Proposition 2.5 from \cite{Zimmer-Lie} with an
additional final claim which we explain.

\begin{propo} \label{K-bundle}
Let $M$ be a manifold with a smooth foliation and a leafwise
Riemannian metric such that every leaf is isometrically covered by
$X_G$. Then there exists a principal right $K$-bundle $\pi_0 : M^*
\rightarrow M$, where $K$ is a maximal compact subgroup of $G$, so
that the $K$-action on $M^*$ extends to a locally free right
$G$-action on $M^*$, which thus defines a foliation whose leaves are
precisely the inverse images under $\pi_0$ of the leaves of $M$.
Furthermore, if $M$ has a dense leaf, then $M^*$ has a dense
$G$-orbit. Finally, any transverse geometric structure defined on
$M$ induces a corresponding $G$-invariant transverse structure on
$M^*$ which is preserved by $\pi_0$.
\end{propo}
\begin{proof}
From \cite{Zimmer-Lie} only the last claim requires some
justification. For that we recall that transverse structures can be
defined by submersions into suitable geometric spaces, as in our
description of Lie structures above. Then, it is enough to compose
with $\pi_0$ the local submersions that define a transverse
structure on $M$. The $G$-invariance follows easily from the first
part of the proposition.\end{proof}

From now on and unless otherwise stated, $M$ will denote a manifold
satisfying the hypotheses of the Main Theorem. Also, $M^*$ will
denote the manifold obtained applying Proposition~\ref{K-bundle} to
$M$. In particular, $M^*$ carries a $G$-invariant transverse Lie
structure modeled in the same group as that for $M$.

We now describe a $G$-invariant pseudoRiemannian metric on $M^*$
that plays a fundamental role to prove our results.

There is a leafwise pseudoRiemannian metric on $M^*$ so that each
leaf is isometrically covered by $G$ with the bi-invariant
pseudoRiemannian metric induced from the Killing form of its Lie
algebra. The latter is easily seen by observing that every leaf of
the foliation on $M^*$ is $G$-equivariantly diffeomorphic to a
quotient $F\backslash G$, for some discrete subgroup $F$.

For $M$ satisfying the hypotheses of the Main Theorem,
Proposition~\ref{K-bundle} shows that $M^*$ carries a $G$-invariant
transverse Lie structure. As we observed in Section~\ref{sec-prel},
this induces a transverse Riemannian structure by choosing a left
invariant Riemannian metric on the simply connected Lie group that
models the transverse Lie structure. In what follows, we fix one
such transverse Riemannian structure.

Our next goal is to show that such leafwise pseudoRiemannian
structure and transverse Riemannian structure can be tied together
into a $G$-invariant pseudoRiemannian metric defined on $M^*$. The
main tools will come from the fact that $\mathcal{F}$ is totally
geodesic. We start by proving that $M^*$ defines a tangential
structure for $M$. From now on, we fix a maximal compact subgroup
$K$ of $G$ and denote $x_0 = eK \in G/K = X_G$.

\begin{propo}\label{M^*-tangential}
If $O_{tg}(\mathcal{F})$ is considered as the bundle of linear
isometries from $T_{x_0}X_G$ onto the fibers of $T\mathcal{F}$, then
the map $M^* \rightarrow O_{tg}(\mathcal{F})$ given by $\varphi
\mapsto d\varphi_{x_0}$ realizes $M^*$ as a tangential $K$-structure
of~$\mathcal{F}$.
\end{propo}
\begin{proof}
Since every isometry is completely determined by its $1$-jet at any
given point, the given map is an embedding. Also, if we consider the
homomorphism from $K$ into $O(T_{x_0}X_G) \simeq O(p)$ ($p =\dim
X_G$) given by $\varphi \mapsto d\varphi_{x_0}$, such embedding is
clearly equivariant and so $M^*$ is a $K$-reduction of
$O_{tg}(\mathcal{F})$. It remains to show that
$\mathcal{H}_\alpha\subset T_\alpha M^*$ for every $\alpha\in M^*$.

Let $\alpha \in M^*$ be given and suppose that $\alpha$ projects to
$x$ in $M$, so that $\alpha : T_{x_0}X_G \rightarrow
T_x\mathcal{F}$. For a given curve $\gamma : [0,1] \rightarrow M$
starting at $x$ and perpendicular to $\mathcal{F}$, let $(\psi_t)_t$
be the elements of horizontal holonomy associated to $\gamma$. Then
for the curve $\overline{\gamma}(t) = d(\psi_t)_x\circ\alpha$ we
have $\overline{\gamma}'(0) \in \mathcal{H}_\alpha$ and every
element in $\mathcal{H}_\alpha$ is of this form; this is a
consequence of the fact that $\mathcal{F}$ is totally geodesic (see
\cite{Blumenthal-Hebda2}). Hence, it suffices to show that
$\overline{\gamma}$ lies in $M^*$ for any such $\alpha$ and
$\gamma$. Let $\varphi\in M^*$ be such that $\varphi(x_0) = x$ and
$d\varphi_{x_0} = \alpha$. In particular, we have
$\overline{\gamma}(t) = d(\psi_t\circ\varphi)_{x_0}$ for every $t$.
We recall that a small enough neighborhood of a locally symmetric
space can be isometrically embedded in a (global) symmetric space.
Hence, if the neighborhoods $V_{\gamma(t)}$ are small enough, for
every $t$ there exists an isometric covering $\varphi_t : X_G
\rightarrow L_{\gamma(t)}$ onto the leaf $L_{\gamma(t)}$ containing
$\gamma(t)$ such that $\varphi_t = \psi_t\circ\varphi$ in a
neighborhood of $x_0$ in $X_G$. In particular, $\overline{\gamma}(t)
= d(\varphi_t)_{x_0}$ and since $\varphi_t \in M^*$ for every $t$,
we conclude that $\overline{\gamma}$ lies in $M^*$ as
required.\end{proof}

From the previous result it follows that the lifted horizontal
bundle $\mathcal{H}$ restricted to $M^*$ is tangent to $M^*$. We
will denote such restriction with $\mathcal{H}^*$. In order to
further understand the properties of $\mathcal{H}^*$ we will need
the following elementary fact about symmetric spaces. For the
definition of standard horizontal vector fields on the frame bundle
of a Riemannian manifold we refer to \cite{KN1}.

\begin{lemma}\label{standard-hor-X_G}
Let $O(X_G)$ be the orthonormal frame bundle of $X_G$ viewed as the
space of linear isometries of $T_{x_0}X_G$ onto the fibers of
$TX_G$. Then, for the embedding $G \hookrightarrow O(X_G)$ given by
$g \mapsto dg_{x_0}$, the standard horizontal vector fields on
$O(X_G)$ are tangent to $G$. Also, if $X = \Gamma\backslash X_G$ is
a Riemannian quotient of $X_G$, where $\Gamma \subset G$ is
discrete, then there is an induced embedding $\Gamma\backslash G
\hookrightarrow O(X)$ so that the standard horizontal vector fields
are tangent to $\Gamma\backslash G$.
\end{lemma}
\begin{proof}
First we observe that the embedding $G \hookrightarrow O(X_G)$ is
$G$-equivariant for the natural left $G$-action on $G$ and the
$G$-action on $O(X_G)$ that lifts from that on $X_G$.

Choose $v\in T_{x_0}X_G$ and denote with $(T_t)_t$ the one-parameter
subgroup of transvections of $X_G$ along the geodesic $\gamma_v$
with initial velocity vector $v$; we recall that a transvection
along a geodesic is an isometry that leaves invariant the geodesic
and whose derivative defines the parallel transport along such
geodesic. Thus $\widehat{\gamma}_v(t) = d(T_t)_{x_0}$ yields the
parallel transport along the geodesic $\gamma_v$ and so it is the
horizontal lift at $\alpha_0$ in $O(X_G)$ with respect to the
Levi-Civita connection, where $\alpha_0 \in O(X_G)$ is the identity
map of $T_{x_0}X_G$. Hence, $\widehat{v} = \widehat{\gamma}_v'(0)$
is horizontal at $\alpha_0$; moreover all horizontal vectors at
$\alpha_0$ are obtained through this construction. Since $T_t \in G$
for every $t$, then $\widehat{\gamma}_v$ lies in $G$ with respect to
the given embedding. Hence, $T_eG$ contains the horizontal subspace
of $O(X_G)$ at $\alpha_0$ defined by the Levi-Civita connection.
Since the embedding $G \hookrightarrow O(X_G)$ is $G$-equivariant
and the $G$-action preserves the Levi-Civita connection on $O(X_G)$
the first part follows. The last claim follows by modding out by the
left $\Gamma$-action and using the properties obtained so
far.\end{proof}

\begin{propo}\label{H^*-invariance}
The lifted horizontal bundle $\mathcal{H}^*$ on $M^*$ is
$G$-invariant.
\end{propo}
\begin{proof}
Let $(X_i)_{i=1}^p$ be the basis for the standard horizontal vector
fields with respect to the Levi-Civita connection along the leaves
of $\mathcal{F}$ and corresponding to the canonical basis of
$\mathbb{R}^p$. Let $(A_i)_{i=1}^l$ be the vertical fields of the
bundle $M^* \rightarrow M$ that come from the right action of the
one-parameter subgroups defined by a basis of $\mathfrak{k}$ (the
Lie algebra of $K$). For such choices, denote $\mathcal{P} = (X_1,
\dots, X_p, A_1,\dots,A_l)$; then, $\mathcal{P}$ consists of
restrictions to $M^*$ of a subset of the natural tangential
parallelism of $O_{tg}(\mathcal{F})$ (see \cite{Cairns1}). Also, by
Proposition~\ref{M^*-tangential}, the set $\mathcal{P}$ defines a
parallelism for the foliation by $G$-orbits in $M^*$.

Also, from the proof of Lemma~\ref{standard-hor-X_G}, the last claim
in Proposition~\ref{K-bundle} and the fact that $G$ preserves the
horizontal standard vector fields on $X_G$, it follows that the
elements in $\mathcal{P}$ are locally given by left invariant vector
fields on $G$. We recall the (elementary) fact that on a Lie group
the left invariant vector fields have flows corresponding to right
actions of one-parameter subgroups. From this fact, again the last
claim of Proposition~\ref{K-bundle} and since $\mathcal{P}$ is a
parallelism, it is easy to see that $\mathcal{P}$ consists of vector
fields whose local flows generate the right $G$-action on $M^*$.

Let $v^* \in \mathcal{H}^*$ be given. Choose curves $\gamma$ in $M$
(perpendicular to $\mathcal{F}$) and $\overline{\gamma}$ in $M^*$ as
in the proof of Proposition~\ref{M^*-tangential}, so that
$\overline{\gamma}'(0) = v^*$ and $\overline{\gamma}(t) =
d(\psi_t)_x\circ\alpha$, where $(\psi_t)_t$ are the elements of
horizontal holonomy associated to $\gamma$, $\alpha$ is the
basepoint of $v^*$ and $x$ is the projection of $\alpha$. Clearly,
the curves $\delta_y(t) = \psi_t(y)$, where $y\in V_x$ with our
previous notation, are the integral curves of a local vector field
$Z^*$ such that $Z^*_\alpha = v^*$; furthermore, $Z^*$ is a section
of $\mathcal{H}^*$ (see \cite{Blumenthal-Hebda2}). Since the local
flow of $Z^*$ is given by maps of the form $y \mapsto \delta_y(t)$
and the latter restrict to the maps $(\psi_t)_t$ along the leaves of
$\mathcal{F}$, it follows that the flow of $Z^*$ fixes the elements
in $\mathcal{P}$. Hence, every element of $\mathcal{P}$ commutes
with $Z^*$, and so the flows given by the elements in $\mathcal{P}$
fix $Z^*$. Since such flows generate the $G$-action, then $G$ maps
$v^*$ into $\mathcal{H}^*$ from which the result follows.\end{proof}

As remarked in the proof of Proposition~\ref{H^*-invariance}, the
elements in $\mathcal{P}$ are locally given by left invariant vector
fields. Hence, by recalling the properties of the metric in $X_G$ in
terms of $G$, we can consider a pseudoRiemannian metric $h_1$ along
the leaves of $\mathcal{F}_{tg}$ in $M^*$ such that the linear span
of the elements in $\mathcal{P}$ are isometric to $\mathfrak{g}$
(the Lie algebra of $G$) with its Killing form. Since the latter is
invariant under the adjoint action of $G$, it is easy to prove that
$h_1$ along the foliation $\mathcal{F}_{tg}$ is $G$-invariant. On
the other hand, from the $G$-invariance of $\mathcal{H}^*$ given by
Proposition~\ref{H^*-invariance}, and since $M$ carries a
$G$-invariant transverse Riemannian structure, there exists a
Riemannian metric $h_2$ on $\mathcal{H}^*$ which is $G$-invariant.
Note that the tangent bundle to the $G$-orbits on $M^*$ is
$T\mathcal{F}_{tg}$. Since $G$ preserves the decomposition $TM^* =
T\mathcal{F}_{tg}\oplus\mathcal{H}^*$ it follows that
$(T\mathcal{F}_{tg},h_1)\oplus(\mathcal{H}^*,h_2)$ defines a
$G$-invariant pseudoRiemannian metric $h$ on $M^*$. From these
remarks we obtain the following result.

\begin{propo} \label{geomstruc1}
Let $M$ satisfy the hypotheses of the Main Theorem. Then, the
manifold $M^*$ carries a $G$-invariant pseudoRiemannian metric $h$
that satisfies the following properties:
\begin{enumerate}
    \item Every $G$-orbit is nondegenerate and the metric induced on
    the $G$-orbits corresponds to the bi-invariant metric on $G$
    given by the Killing form of its Lie algebra.
    \item The metric $h$ defines transverse
    Riemannian structure that corresponds to the one given by the
    transverse Lie structure as described in Section \ref{sec-prel}.
\end{enumerate}
\end{propo}

In what follows, $h$ will denote the pseudoRiemannian metric on
$M^*$ described by the previous proposition.

\section{Proof of the Main Theorem}
\label{sec-main-proofs} Let $h$ be the $G$-invariant
pseudoRiemannian metric on $M^*$ constructed in Section
\ref{sec-M^*}. Following \cite{Quiroga-Annals} we define $n_0$ and
$m_0$ as the maximal dimensions of the null tangent subspaces for
$G$ and $M^*$, respectively; for $G$ we consider the bi-invariant
pseudoRiemannian metric given by the Killing form of its Lie
algebra. Since the metric $h$ restricted to the $G$-orbits is
precisely this bi-invariant metric and $h$ is Riemannian in the
orthogonal complement to the tangent bundle to the $G$-orbits, it
follows easily that $n_0 = m_0$. On the other hand, by Proposition
\ref{K-bundle} there is a dense $G$-orbit in $M^*$. Hence the
hypothesis of Theorem A in \cite{Quiroga-Annals} are satisfied. It
follows from this reference the existence of a finite covering
$\pi^* : \widehat{M^*} \rightarrow M^*$ over which the $G$-action
can be lifted so that there is a $G$-equivariant diffeomorphism:
$$
    \varphi^* : \Gamma\backslash(N \times G) \rightarrow
    \widehat{M^*},
$$
where $N$ is a homogeneous Riemannian manifold and $\Gamma \subset
\operatorname{Iso}(N)\times G$. That $N$ is Riemannian follows from
the fact that the orthogonal complement to the $G$-orbits is
Riemannian in our construction of the metric $h$ on $M^*$. Note that
without loss of generality we can assume that $N$ is simply
connected by replacing the covering $N \times G \rightarrow
\widehat{M^*}$ with the its composition with the universal covering
map of $N$ in the first factor.

From this discussion, the manifolds $N \times \widetilde{G}$ and
$\widetilde{M^*}$ are $\widetilde{G}$-equivariantly diffeomorphic.
Let $H$ be a Lie group that models the transverse Lie structure of
$\widehat{M^*}$. By Theorem \ref{Zimmer-hol-thm} the group $H$ is
semisimple. Choose a development $D : N \times \widetilde{G} \cong
\widetilde{M^*} \rightarrow H$. In other words, $D$ is a submersion
equivariant with respect to the holonomy representation $\rho :
\pi_1(\widehat{M^*}) \rightarrow H$. Clearly, we can assume that $H$
is simply connected.

The foliation in $\widehat{M^*}$ is given by the $G$-orbits, and so
the foliation in $N \times \widetilde{G} \cong \widetilde{M^*}$ is
given by the $\widetilde{G}$-orbits. Since $D$ is a development, it
follows that the connected components of the fibers of $D$ are the
$\widetilde{G}$-orbits. Hence, the restriction $D|_{N\times\{e\}}$
defines a local diffeomorphism $N \rightarrow H$. Moreover, by the
choice of $N$ and its Riemannian metric as given in
Proposition~\ref{geomstruc1}, such map is a local isometry, hence
(by Corollary 29 in page 202 in \cite{ONeill}) a covering map and so
a diffeomorphism. It follows that $N$ admits a Lie group structure
so that the map $D|_{N\times\{e\}} : N \rightarrow H$ is an
isomorphism. Hence, if we consider $N$ endowed with such Lie group
structure, then the canonical projection $N \times \widetilde{G}
\rightarrow N$ is a development for the transverse Lie structure of
$\widehat{M^*}$. In particular, we can replace the Riemannian
manifold $N$ with the Lie group $H$ carrying a suitable left
invariant Riemannian structure.

At the same time, since the deck transformation group of the
covering $\widetilde{G} \rightarrow G$ is a central subgroup of
$\widetilde{G}$ and $\pi_1(\widehat{M^*}) \subset
\operatorname{Iso}(H) \times \widetilde{G}$, then using the
diffeomorphism $\varphi^*$ above it follows that the natural
projection $H\times G \rightarrow H$ defines a development for the
transverse Lie structure on $\widehat{M^*}$ as well.

On the other hand, if $\gamma \in \Gamma$, then we can write $\gamma
= (\gamma_1, \gamma_2)$, where $\gamma_1$ is a diffeomorphism of $H$
and $\gamma_2\in G$. Since $\pi$ is equivariant with respect to the
corresponding holonomy representation $\rho : \Gamma \rightarrow H$,
we conclude that:
$$
    \gamma_1(x)  = D(\gamma_1 x, \gamma_2 g) = D(\gamma(x,g)) =
            \rho(\gamma) D(x,g) = y x
$$
is satisfied for some $y \in H$, depending only on $\gamma$, and
every $(x,g)\in H \times G$. It follows that $\gamma_1$ is the left
translation by an element in $H$ and with respect to this we have
$\gamma_1 \in H$, thus showing that $\Gamma$ is a subgroup of
$H\times G$.

Next, if we consider the finite covering:
$$
       \pi^* : \Gamma\backslash(H\times G) \rightarrow M^*
$$
just constructed, then for $K$ a maximal compact subgroup of $G$ as
before, we observe that taking the quotient on this map by the
$K$-action on the right, we obtain a finite covering map:
$$
       \pi : \Gamma\backslash(H\times G/K) \rightarrow M,
$$
where $G/K$ is the symmetric space that isometrically covers the
leaves of the foliation in $M$. Hence, from the previous discussion
most of the claims of the Main Theorem have now been proved. It
remains to show that we can choose our covering $\pi$ so that $H$ is
replaced by a finite center quotient group; this is achieved by a
modding out by a suitable central subgroup of $H$ through a standard
procedure which already appears in Section 6 of
\cite{Quiroga-Annals}. Finally, the arithmeticity of $\Gamma$ is a
consequence of Margulis' arithmeticity theorem.

\bibliographystyle{amsplain}

\end{document}